
\input amstex
\documentstyle{amsppt}
\magnification=1200
\define \R{\Bbb R}
\define \ir{\int\limits_{\Bbb R}}
\define \ft{Fourier transform}
\define \fs{Fourier series}
\define \T{\Bbb T}
\define \iT{\int\limits_{\Bbb T}}
\define \hf{\hat f}
\define \hg{\hat g}
\define \la{\lambda}
\define \La{\Lambda}
\define \cb{\Cal B}
\define \ch{\Cal H}
\define \ca{\Cal A}

\topmatter
\author E. Liflyand \endauthor
\title A space of multipliers on L \endtitle
\address Department of Mathematics and Computer Sciences, Bar-Ilan University,
52900 Ramat-Gan, Israel   \endaddress
\email liflyand\@bimacs.cs.biu.ac.il \endemail
\abstract Conditions for a function (number sequence) to be a
multiplier on the space of integrable functions on $\Bbb R$ ($\Bbb T$)
are given. This generalizes recent results of Giang and Moricz. \endabstract
\endtopmatter
\document
\vskip 3cm
\centerline{\bf BIMACS - 9503}
\vskip 3cm
\centerline{\bf Bar-Ilan University - 1995}
\vfill\eject
\head 1. Introduction \endhead
In their recent paper \cite{GM}, D\u ang V\~u Giang and F. M\'oricz
gave a family of spaces, so that each element of such space is a
multiplier on $$L=\{f:||f||_L =\ir |f(x)|\,dx < \infty\},$$
where $\R = (-\infty,\infty).$ These results are obtained both
in periodic and non-periodic cases. Proofs are strongly based
on some sufficient conditions for a function to have an integrable
\ft , or for a trigonometric series to be a \fs\ of an integrable
function, respectively.

More general conditions of such type were given in our recent paper
\cite{L}. Thus we can give less restrictive multiplier conditions,
that is the results of \cite{GM} follow from ours immediately.
Some of our notation is the same as in \cite{GM}. It allows us to
compare easily our results.

Let $f$ be an integrable function on $\R$, and
$$\hf (x) = \ir f(t) e^{-ixt}\, dt $$ be its \ft .

We say that a measurable bounded function $\la$ is an $M$-multiplier
if for every $f\in L$ there exists a function $g\in L$ such that
$$\la (t) \hf (t)= \hg (t) .  \tag1$$
The norm of the corresponding operator $\La : L\to L$ which assigns
to each $f\in L$ the function $\La f=g$ accordingly to (1) may be
calculated as usually:
$$||\La ||_M = \sup_{||f||_L \le 1} {||\La f||_L \over ||f||_L}.$$

One may consider the space $\cb$ of absolutely continuous functions
on $\R$, bounded over $\R$, and endowed with the norm
$$||\la||_{\cb}=||\la||_B +\ir |\la '(t)|\,dt ,$$
where $||\la ||_B =\sup_{t\in \R}|\la (t)| .$ There exist functions in $\cb$
which are not multipliers (see, e.g.,  \cite{T}, p.170-172). Thus
it is interesting to study some subspaces of $\cb$ which are the
spaces of multipliers.

\head 2. Description of the space of multipliers \endhead
We introduce a subspace of $\cb$ denoted by $\ch$ and defined as follows.
Let $$S_f =\int_0^{\infty} \biggl|\int\limits_{|t|\le {u\over 2}}
{f(u-t) - f(u+t) \over t}\,dt \biggr|\,du .$$
Then, denoting $$\ca =\int\limits_0^{\infty} {|\la (t)-\la (-t)|
\over t}\, dt ,$$ we set
$$\ch =\{\la : ||\la||_{\ch} =||\la||_{\cb} +S_{\la'} +\ca  <\infty \}.$$

Our main result is the following
\proclaim{Theorem 1} If $\la \in \ch$ then $\la$ is an $M$-
multiplier, and $$||\La||_M \le C||\la||_{\ch} . \tag2$$ \endproclaim

Here and in what follows $C$ will mean absolute constants, and $C$
with indices, say $C_p ,$ will denote some constants depending only
on the indices mentioned. The same letter may denote constants
different in different places.

Let us make some remarks on the space $\ch$ and compare Theorem 1
with earlier results.

Let $f_+$ be the odd continuation of the part of a function $f$
supported on $[0,\infty),$ and $f_-$ be the odd continuation of
the part of $f$ supported on $(-\infty,0].$

Let, further, $ReH$ be the Hardy space with the norm
$$||f||_H = ||f||_L + ||\tilde f||_L <\infty ,$$
where $$\tilde f (x)={1\over \pi}\ir {f(t)\over x-t}\,dt $$
is the Hilbert transform of $f$.

We say that $f\in H_s$ if $f_+ \in ReH$ and $f_- \in ReH$.
It was proved in \cite{L} that $||f||_1 + S_{f} <\infty$
is equivalent to the fact that $f\in H_s .$

In \cite{GM}, the main theorem is similar to our Theorem 1, but with
one of the spaces of a family $\{\cb_p :1\le p<\infty\}$ instead of
$\ch .$ This family was defined as follows. For $1<p<\infty$ set
$$\ca_q f= \int\limits_0^{\infty} \biggl({1\over u} \int\limits_{u\le |t|
\le 2u} |f(t)|^q \,dt \biggr)^{{1\over q}}\,du $$
where ${1\over p}+{1\over q}=1$ here and in what follows, while for $p=1$
$$\ca_{\infty} f =\int\limits_0^{\infty} \operatornamewithlimits
{ess\,sup}\limits_{u\le |t| \le 2u} |f(t)|\,du .$$
Then for $1\le p<\infty$
$$\cb_p =\{\la : ||\la||_{\cb_p} =||\la||_B +\ca_q \la' +\ca <\infty\}.$$

A theorem claimed as the main one in \cite{GM} may be formulated
like Theorem 1, with $\cb_p$ instead of $\ch .$ Its proof, after
some simple computations, follows from the following
\proclaim{Lemma} {\rm (See Lemma 1 in \cite{GM}.)} If a function $\la$
is locally absolutely continuous, satisfies the condition
$$\lim_{|t|\to\infty} \la (t) = 0 , \tag3$$
and for some $1<q\le\infty$ we have $\ca_q <\infty ,$ then $\hat\la$
belongs to $L$ if and only if $\ca <\infty .$ Furthermore,
$$||\hat \la||_L \le \ca +C_q \ca_q .$$          \endproclaim

Indeed,
the fact that $\la'\in L^1$ yields easily that $\la$ has finite limits
$l_+$ and $l_-$ at $+\infty$ and $-\infty,$ respectively. When they
coincide, $l_+=l_-=l$ (and this follows from $\ca<\infty$), one can take
$$\la_0(t)=\la(t)-l$$ and $$\la_1(t)=\hat\la_0(-t).$$
If $\la_1\in L$ then $\hat\la_1=\la_0,$ and taking
$$\Lambda f=lf+\la_1*f$$  one gets
$$\align (\Lambda f)(t) &=l\hat f(t)+\hat\la_1\hat f(t)\\
&=l\hat f(t)+\la_0(t)\hat f(t)=\la(t)\hat f(t), \endalign$$
and $\la$ is a multiplier on $L.$

Theorem 1 may be proved analogously, with application of one our result
from \cite{L} instead of Lemma.
\proclaim{Theorem A} {\rm (See \cite{L}, Theorem 2.)} Let $\la$ be a locally
absolutely continuous function, satisfying {\rm (3)}. Then for $|x|>0$ we have
$$\hat\la (x)={i\over x}\biggl(\la ({\pi\over 2|x|})-\la (-{\pi\over 2|x|})
\biggr) + \theta\gamma (x) ,$$ where $|\theta|\le C ,$ and
$$\ir |\gamma (x)|\,dx \le ||\la'||_L + S_{\la'} .$$   \endproclaim

It is obvious that, in conditions of Theorem A, $\hat\la\in L$ iff
$\ca <\infty .$ As it was said above, in order to prove Theorem 1
it remains to repeat the proof of Theorem 1 from \cite{GM} using
Theorem A instead of Lemma 1.

Indeed, the following embeddings are almost obvious:
$$\cb_1\subset\cb_{p_1}\subset\cb_{p_2}\subset\cb ,\qquad 1<p_1<p_2<\infty,$$
while the following fact, proved in \cite{L}, is not so clear:
$$\cb_p\subset\ch ,\qquad 1\le p <\infty .$$
Therefore, the main result in \cite{GM} is contained in Theorem 1.

That  $\la\in\ca_q$ may not be in $ReH$
can be seen from the following counterexample. Let $\la (x)={1\over 1+x^2}.$
We have $\ca_q \la<\infty.$
Nevertheless $\tilde{\la}(x) ={x\over 1+x^2}$ (see, e.g., \cite{BN}, p.519),
and $\tilde{\la}\not\in L.$

Easy sufficient condition for an even function $\la$ defined on $[0,\infty)$
to be a Fourier multiplier, due to \cite{BN}, p.248, has the following
relation (so-called quasiconvexity)
$$\int\limits_0^{\infty} t|d\la'(t)|\,<\infty$$
as the main part. It is easy to verify that this condition is more
restrictive than $S_{\la'}<\infty .$ Indeed, we have
$$\align & S_{\la'}=\int\limits_0^{\infty}\biggl|
\int\limits_0^{{u\over 2}}\, {dx\over x}\int\limits_{u-x}^{u+x}
\, d\la'(t)\biggr|\,du \\ & \le \int\limits_0^{\infty} \,du
\int\limits_{{u\over 2}}^{{3u\over 2}}\,|d\la'(t)|\ln{u\over 2|u-t|} =
\ln 3 \int\limits_0^{\infty} t|d\la'(t)| . \endalign$$

\head 3. The case of Fourier series \endhead
Analogous results for the case of \fs\ were obtained in \cite{GM}
as well. We can generalize these results in the same manner as in the
case of \ft s.

Let now $L$ be the space of all complex-valued $2\pi$-periodic
functions integrable over $\T =(-\pi,\pi] ,$ and
$$||f||_L =\iT |f(x)|\,dx .$$      Let
$$\hf (k) ={1\over 2\pi}\iT f(x)e^{-ikx}\,dx, \qquad k=0,\pm 1,\pm 2,...$$
be the Fourier coefficients of the function $f.$

A bounded sequence $\{\la =\la (k)\}$, with $||\la||_m =\sup |\la (k)|
<\infty ,$ is called an $M$-multiplier if for every $f\in L$ there
exists a function $g\in L$ such that
$$\la (k) \hf (k)=\hg (k),\qquad k=0,\pm 1,\pm 2,... \tag4$$
As above (4) assigns a bounded linear operator $\La ,$ and it is
worth studying spaces of multipliers which are subspaces of the space
$$bv =\{\la :||\la||_{bv} =||\la||_m +||\Delta\la||_1 <\infty\} ,$$
where $||\cdot||_1$ is the norm in $l^1 ,$ and
$$\Delta\la (k) =\cases \la (k)-\la (k+1) & \text{if} \quad k\ge 0 , \\
\la (k)-\la (k-1) & \text{if} \quad k<0 . \endcases$$
Again a sequence $\la\in bv$ exists which is not a multiplier
(see \cite{Z}, Vol.1, p.184).

We introduce a subspace of $bv$
$$h=\{\la: ||\la||_h =||\la||_{bv} +s_{\la} + a <\infty\}$$
where $$s_{\la}=\sum\limits_{m=2}^{\infty}\,\biggl|\,\sum\limits_{k=1}^
{[{m\over 2}]} {\Delta\la (m-k) -\Delta\la (m+k) \over k}\biggr|$$
and $$a=\sum\limits_{k=1}^{\infty} {|\la (k)-\la (-k)|\over k} .$$
Note that the condition $s_{\la} <\infty$ is called the
Boas-Telyakovskii condition (see, e.g., \cite{T1}).

The analog of Theorem 1 for \fs\ may be formulated as follows.
\proclaim{Theorem 2} If $\la\in h$ then $\la$ is an $M$-multiplier,
and $$||\La||_M \le C||\la||_h .$$     \endproclaim

The proof again may be reduced to the proof of the multiplier
theorem for \fs\ in \cite{GM} with application of the following corollary
to Theorem A instead of corresponding weaker result in \cite{GM}
(see Lemma 3).

Let $\ell (x) =\la(k) +(k-x)\Delta \la(k)$ for $x\in[k-1,k]$, with
$\lim\limits_{|k|\to\infty} \la(k) =0.$
\proclaim{Theorem B} {\rm (see \cite{L}, Theorem 5).} For every $y,$
$0<|y|\le\pi ,$ $$\sum\limits_{k=-\infty}^{\infty} \la(k) e^{iky} =
{i\over y} \biggl(\ell({\pi\over 2|y|})-\ell(-{\pi\over 2|y|})\biggr)
+\theta\gamma (y) \tag4$$ where $\theta \le C,$ and
$$\iT |\gamma (y)|\,dy \le ||\Delta\la||_1 + s_{\la} .$$
\endproclaim

This is a somewhat stronger form of Telyakovskii's result in [T1].

It is obvious now that the function $\ell$, having the sequence
$\la$ as its Fourier coefficients, is integrable over $\T$ when
$\la\in h .$ Thus it is enough to substitute this result
for Lemma 3 in \cite{GM}, and so changed proof establishes
Theorem 2.

Analogously to the case of \ft s, the multiplier properties
in the case of \fs\ were proved in [GM] for a family of sequences
$\{bv_p , 1\le p<\infty\}.$ Each family is a subspace of $bv$
and is defined as follows. Let $I_n =\{2^n,2^n+1,...,2^{n+1}-1\},$ and
$$a_q =\sum\limits_{n=0}^{\infty}2^n\biggl(2^{-n}\sum\limits_
{|k|\in I_n} |\Delta\la(k)|^q\biggr)^{1\over q} $$
for $1<p<\infty,$ while for $p=1$
$$a_{\infty} =\sum\limits_{n=0}^{\infty} \max_{|k|\in I_n} |\Delta\la(k)|.$$
Then for $1\le p<\infty$
$$bv_p =\{\la:||\la||_{bv_p} =||\la||_m +a_q +a <\infty \}.$$
It is well known that
$$bv_1\subset bv_{p_1}\subset bv_{p_2}\subset bv,\qquad 1<p_1<p_2<\infty$$
but for us more important is that for all $p$ $$bv_p \subset h.$$
This was proved for $p=1$ by Telyakovskii \cite{T2}, and for $p>1$
by Fomin \cite{F}. Therefore the result for \fs\ in \cite{GM}
is contained in Theorem 2 as the partial case.

\enddocument